\documentclass{article}
\usepackage{amssymb}
\usepackage{amsmath}
\usepackage{amsfonts}
\usepackage{hyperref}

\newtheorem{theorem}{Theorem}

\newtheorem{remark}[theorem]{Remark}

\input{tcilatex}
\numberwithin{theorem}{section}
\numberwithin{equation}{section}

\begin{document}

\title{Riemann-Lagrange Geometric Dynamics for the Multi-Time Magnetized
Non-Viscous Plasma}
\author{Mircea Neagu}
\date{}
\maketitle

\begin{abstract}
In this paper, using Riemann-Lagrange geometrical methods, we construct a
geometrical model on 1-jet spaces for the study of multi-time relativistic
magnetized non-viscous plasma, characterized by a given
energy-stress-momentum distinguished (d-) tensor. In that arena, we give the
conservation laws and the continuity equations for multi-time plasma. The
partial differential equations of the stream sheets (the equivalent of
stream lines in the classical semi-Riemannian geometrical approach of
plasma) for multi-time plasma are also written.
\end{abstract}



\section{Introduction}

\hspace{5mm}During that so-called the radiation epoch, in which photons are
strongly coupled with the matter, the interactions between the various
constituents of the Universal matter include radiation-plasma coupling,
which is described by the plasma dynamics. Although it is not traditional to
characterize the radiation epoch by the dominance of plasma interactions,
however, it may be also called the plasma epoch (please see \cite{Kleidis}).
This is because the electromagnetic interaction dominates all the four
fundamental physical forces (electrical, magnetic, gravitational and
nuclear).

In the present days, the Plasma Physics is an well established field of
Theoretical Physics, although the formulation of magnetohydrodynamics in a
curved space-time is a relatively new development (please see Punsly \cite%
{Punsly}). The MHD processes in an isotropic space-time are intensively
studied by a lot of physicists. For example, the MHD equations in an
expanding Universe are investigated by Kleidis, Kuiroukidis, D. Papadopoulos
and Vlahos in \cite{Kleidis}. Con\-si\-dering the interaction of the
gravitational waves with the plasma in the presence of a weak magnetic
field, D.B. Papadopoulos also investigates the relativistic hydromagnetic
equations \cite{Papadopoulos}. The electromagnetic-gravitational dynamics
into plasmas with pressure and viscosity is studied by Das, DeBenedictis,
Kloster and Tariq in the paper \cite{Das}. In their paper, the authors
derive the relativistic Navier-Stokes equations that govern plasma.

It is important to note that all preceding physical studies are done on an
isotropic four-dimensional space-time, represented by a semi- (pseudo-)
Riemannian space with the signature $(+,+,+,-)$. Consequently, the
Riemannian ge\-o\-me\-tri\-cal methods are used as a pattern over there.

Geometrically speaking, using the Finlerian geometrical methods, the plasma
dynamics was extended on non-isotropic space-times by V. G\^{\i}r\c{t}u and
Ciu\-bo\-ta\-ri\-u in the paper \cite{Girtu-Ciubotariu}. More general, after
the development of Lagrangian geometry on tangent bundle, due to Miron and
Anastasiei \cite{Mir-An}, the generalized Lagrange geometrical objects
describing the relativistic magnetized plasma were studied by M. G\^{\i}r%
\c{t}u, V. G\^{\i}r\c{t}u and Postolache in the paper \cite%
{Girtu-Girtu-Posto}.

According to Olver's opinion \cite{Olver}, we appreciate that the $1$-jet
fibre bundle is a basic object in the study of classical and quantum field
theories. For such a reason, using as a pattern the Miron-Anastasiei's
geometrical ideas \cite{Mir-An}, the author of this paper recently developed
in the paper \cite{Neagu-Forum} that so-called the "multi-time
Riemann-Lagrange geometry"\textit{\ }on 1-jet spaces, in the sense of
d-connections, d-torsions, d-curvatures, gravitational and electromagnetic
geometrical theories. We would like to point out that the geometrical
construction on $1$-jet spaces exposed in the article \cite{Neagu-Forum} was
initiated by Asanov in \cite{Asanov} and further developed by the author of
this paper. In this geometrical context, the aim of this paper is to create
a multi-time extension on $1$-jet spaces of the geometrical objects that
characterize plasma in semi-Riemannian and Lagrangian approaches. Thus, we
introduce the energy-stress-momentum d-tensor of the "multi-time plasma" and
we give the geometrical-physical equations which govern it.

\section{The semi-Riemannian geometrical approach. Plas\-ma in isotropic
space-times}

\hspace{5mm}Let $\mathcal{SR}_{n}=\left( M^{n},\varphi _{ij}(x)\right) $ be
a semi-Riemannian manifold, where $M^{n}$ is an $n$-dimensional smooth
manifold, whose coordinates are $x=(x^{i})_{i=\overline{1,n}}$, and $\varphi
_{ij}(x)$ is a semi-Riemannian metric having a constant signature. From a
physical point of view, $\varphi _{ij}(x)$ play the role of gravitational
potentials. Note that, throughout this paper, the latin letters run from $1$
to $n$ and the Einstein convention of summation is assumed. Let us consider
the Christoffel symbols of the semi-Riemannian metric $\varphi _{ij}$, which
are given by%
\begin{equation}
\gamma _{jk}^{i}=\frac{\varphi ^{im}}{2}\left( \frac{\partial \varphi _{jm}}{%
\partial x^{k}}+\frac{\partial \varphi _{km}}{\partial x^{j}}-\frac{\partial
\varphi _{jk}}{\partial x^{m}}\right) .  \label{Christoffel-Riemann}
\end{equation}%
The Christoffel symbols (\ref{Christoffel-Riemann}) produce the Levi-Civita
covariant de\-ri\-va\-tive%
\begin{equation*}
T_{kl...;p}^{ij...}=\frac{\partial T_{kl...}^{ij...}}{\partial x^{p}}%
+T_{kl...}^{mj...}\gamma _{mp}^{i}+T_{kl...}^{im...}\gamma
_{mp}^{j}+...-T_{ml...}^{ij...}\gamma _{kp}^{m}-T_{km...}^{ij...}\gamma
_{lp}^{m}-...,
\end{equation*}%
where%
\begin{equation*}
T=T_{kl...}^{ij...}(x)\frac{\partial }{\partial x^{i}}\otimes \frac{\partial 
}{\partial x^{j}}\otimes dx^{k}\otimes dx^{l}\otimes ...
\end{equation*}%
is an arbitrary tensor on $M$.

\begin{remark}
The Levi-Civita covariant derivative has the metrical properties%
\begin{equation*}
\varphi _{ij;k}=\varphi _{\text{ \ };k}^{ij}=0.
\end{equation*}
\end{remark}

To define the energy-stress-momentum tensor $\mathcal{T}=\mathcal{T}%
_{ij}(x)dx^{i}\otimes dx^{j}$, that characterize the relativistic magnetized
non-viscous plasma, we need the following geometrical objects \cite{Das}, 
\cite{Girtu-Girtu-Posto}:

\begin{enumerate}
\item the unit velocity-field of a test particle, given by%
\begin{equation*}
U=u^{i}(x)\frac{\partial }{\partial x^{i}},
\end{equation*}%
where, if we denote $u_{i}=\varphi _{im}u^{m}$, then we have $u_{i}u^{i}=1$.
Note that, physically speaking, if $V=v^{i}(x)\left( \partial /\partial
x^{i}\right) $ is the fluid's space-like velocity vector, then we have%
\begin{equation*}
u^{i}=\frac{v^{i}}{\sqrt{\varphi _{rs}v^{r}v^{s}}};
\end{equation*}

\item the $2$-form of the (electric field)-(magnetic induction) is given by%
\begin{equation*}
H=H_{ij}(x)dx^{i}\wedge dx^{j}
\end{equation*}%
and the $2$-form of the (electric induction)-(magnetic field) is given by%
\begin{equation*}
G=G_{ij}(x)dx^{i}\wedge dx^{j}.
\end{equation*}%
Note that, in physical applications, one takes $H=-G=F/\sqrt{\mu _{0}}$,
where $F$ is the electromagnetic field and $\mu _{0}$ is the electromagnetic
permeability constant.

\item the Minkowski energy tensor of the electromagnetic field inside the
plasma is given by the tensor $E=E_{ij}(x)dx^{i}\otimes dx^{j}$, whose
components are%
\begin{equation*}
E_{ij}=\frac{1}{4}\varphi _{ij}H_{rs}G^{rs}+\varphi ^{rs}H_{ir}G_{js},
\end{equation*}%
where $G^{rs}=\varphi ^{rp}\varphi ^{sq}G_{pq}.$ In order to obey the
relativistic Lorentz equation of motion for a charged test particle, the
following Lorentz condition is required \cite{Das}:%
\begin{equation}
E_{i;m}^{m}u^{i}=0,  \label{Lorentz-Conditon-Riemann}
\end{equation}%
where $E_{i}^{m}=\varphi ^{mp}E_{pi}$. Obviously, using the notations $%
H_{r}^{m}=\varphi ^{mp}H_{pr}$ and $G_{i}^{r}=\varphi ^{rs}G_{si}$, then we
have%
\begin{equation*}
E_{i}^{m}=\frac{1}{4}\delta _{i}^{m}H_{rs}G^{rs}-H_{r}^{m}G_{i}^{r},
\end{equation*}%
where $\delta _{i}^{m}$ is the Kronecker symbol.
\end{enumerate}

In this physical context, the components of the energy-stress-momentum
tensor of plasma are defined by (please see \cite{Das}, \cite%
{Girtu-Ciubotariu}, \cite{Girtu-Girtu-Posto})%
\begin{equation}
\mathcal{T}_{ij}=\left( \mathbf{\rho }+\frac{\mathbf{p}}{c^{2}}\right)
u_{i}u_{j}+\mathbf{p}\varphi _{ij}+E_{ij},  \label{stress-Riemann}
\end{equation}%
where $c=$ const. is the speed of light, $\mathbf{p}=\mathbf{p}(x)$ is the
hydrostatic pressure and $\mathbf{\rho }=\mathbf{\rho }(x)$ is the proper
mass density of plasma.

In the Riemannian framework of plasma, it is postulated that the following 
\textit{conservation laws} for the components (\ref{stress-Riemann}) are
true:%
\begin{equation}
\mathcal{T}_{i;m}^{m}=0,  \label{Cons-Riemann}
\end{equation}%
where%
\begin{equation*}
\mathcal{T}_{i}^{m}=\varphi ^{mp}\mathcal{T}_{pi}=\left( \mathbf{\rho }+%
\frac{\mathbf{p}}{c^{2}}\right) u^{m}u_{i}+\mathbf{p}\delta
_{i}^{m}+E_{i}^{m}.
\end{equation*}%
By direct computations, the conservation equations (\ref{Cons-Riemann})
become%
\begin{equation}
\left[ \left( \mathbf{\rho }+\frac{\mathbf{p}}{c^{2}}\right) u^{m}\right]
_{;m}u_{i}+\left( \mathbf{\rho }+\frac{\mathbf{p}}{c^{2}}\right)
u^{m}u_{i;m}+\mathbf{p}_{,i}-\varphi _{ir}\mathcal{F}^{r}=0,
\label{cons-local-Riemann}
\end{equation}%
where $\mathbf{p}_{,i}=\partial \mathbf{p}/\partial x^{i}$ and $\mathcal{F}%
^{r}=-\varphi ^{rs}E_{s;m}^{m}$ is the \textit{Lorentz force}.

Contracting the conservation equations (\ref{cons-local-Riemann}) with $%
u^{i} $ and taking into account the Lorentz condition (\ref%
{Lorentz-Conditon-Riemann}), we find the \textit{continuity equation} of
plasma, namely%
\begin{equation}
\left[ \left( \mathbf{\rho }+\frac{\mathbf{p}}{c^{2}}\right) u^{m}\right]
_{;m}+\mathbf{p}_{,m}u^{m}=0,  \label{Continuity-Riemann}
\end{equation}%
where we also used the equalities%
\begin{equation*}
0=u_{i}u_{;m}^{i}=\frac{1}{2}\left( u_{i}u^{i}\right) _{,m}=-u_{i;m}u^{i},
\end{equation*}%
the comma symbol "$_{,m}$" representing the partial derivative $\partial
/\partial x^{m}$.

Replacing the continuity law (\ref{Continuity-Riemann}) into the
conservation equations (\ref{cons-local-Riemann}), we find the \textit{%
relativistic Euler equations} for plasma, namely%
\begin{equation}
\left( \mathbf{\rho }+\frac{\mathbf{p}}{c^{2}}\right) u_{i;m}u^{m}-\mathbf{p}%
_{,m}\left( u^{m}u_{i}-\delta _{i}^{m}\right) -\varphi _{im}\mathcal{F}%
^{m}=0.  \label{Euler-Riemann}
\end{equation}

If we put now $u^{m}=dx^{m}/ds$ into Euler equations (\ref{Euler-Riemann}),
we find out the equations of the \textit{stream lines} of plasma, which are
given by the second order DE system%
\begin{equation*}
\frac{d^{2}x^{k}}{ds^{2}}+\left[ \gamma _{rm}^{k}-\frac{c^{2}}{\mathbf{p}+%
\mathbf{\rho }c^{2}}\delta _{r}^{k}\mathbf{p}_{,m}\right] \frac{dx^{r}}{ds}%
\frac{dx^{m}}{ds}=\frac{c^{2}}{\mathbf{p}+\mathbf{\rho }c^{2}}\left[ 
\mathcal{F}^{k}-\varphi ^{km}\mathbf{p}_{,m}\right] ,
\end{equation*}%
where $s$ is the natural parameter of the smooth curve $c=\left(
x^{k}(s)\right) _{k=\overline{1,n}}.$

\section{Generalized Lagrangian geometrical a\-pproach. Plasma in
non-isotropic space-times}

\hspace{5mm}Let $\mathcal{GL}^{n}=\left( M^{n},\text{ }g_{ij}(x,y),\text{ }%
N_{j}^{i}(x,y)\right) $ be a generalized Lagrange space (for more details,
please see Miron and Anastasiei \cite{Mir-An}). Let us consider that the
tangent bundle $TM$, as smooth manifold of dimension $2n$, has the local
coordinates $(x^{i},y^{i})_{i=\overline{1,n}}$. Then, $g_{ij}(x,y)$ is a
metrical d-tensor on $TM$, which is symmetrical, non-degenerate and has a
constant signature on $TM\backslash \{0\}$. The local coefficients $%
N_{j}^{i}(x,y)$ are the components of a nonlinear connection $N$ on $TM$.
The nonlinear connection $N=\left( N_{j}^{i}\right) $ produces on $TM$ the
following dual adapted bases of d-vectors and d-covectors:%
\begin{equation*}
\left\{ \frac{\delta }{\delta x^{i}},\frac{\partial }{\partial y^{i}}%
\right\} \subset \mathcal{X}(TM),\qquad \left\{ dx^{i},\delta y^{i}\right\}
\subset \mathcal{X}^{\ast }(TM),
\end{equation*}%
where%
\begin{equation*}
\frac{\delta }{\delta x^{i}}=\frac{\partial }{\partial x^{i}}-N_{i}^{m}\frac{%
\partial }{\partial y^{m}},\qquad \delta y^{i}=dy^{i}+N_{m}^{i}dx^{m}.
\end{equation*}%
Note that the d-tensors on the tangent bundle $TM$ behave like classical
tensors. For example, on $TM$ we have the global metrical d-tensor%
\begin{equation*}
\mathbb{G}=g_{ij}dx^{i}\otimes dx^{j}+g_{ij}\delta y^{i}\otimes \delta y^{j},
\end{equation*}%
which is usually endowed with the physical meaning of non-isotropic
gravitational potential.

Following the geometrical ideas of Miron and Anastasiei from \cite{Mir-An},
the ge\-ne\-ra\-lized Lagrange space $\mathcal{GL}^{n}$ produces the Cartan
canonical $N$-linear connection%
\begin{equation*}
C\Gamma (N)=\left( L_{jk}^{i},\text{ }C_{jk}^{i}\right) ,
\end{equation*}%
where%
\begin{equation}
\begin{array}{l}
\medskip L_{jk}^{i}=\dfrac{g^{im}}{2}\left( \dfrac{\delta g_{jm}}{\delta
x^{k}}+\dfrac{\delta g_{km}}{\delta x^{j}}-\dfrac{\delta g_{jk}}{\delta x^{m}%
}\right) , \\ 
C_{jk}^{i}=\dfrac{g^{im}}{2}\left( \dfrac{\partial g_{jm}}{\partial y^{k}}+%
\dfrac{\partial g_{km}}{\partial y^{j}}-\dfrac{\partial g_{jk}}{\partial
y^{m}}\right) .%
\end{array}
\label{Cartan-Lagrange}
\end{equation}%
Further, the Cartan linear connection $C\Gamma (N)$, given by (\ref%
{Cartan-Lagrange}), induces the horizontal ($h-$) covariant derivative%
\begin{equation*}
D_{kl...|p}^{ij...}=\frac{\delta D_{kl...}^{ij...}}{\delta x^{p}}%
+D_{kl...}^{mj...}L_{mp}^{i}+D_{kl...}^{im...}L_{mp}^{j}+...-D_{ml...}^{ij...}L_{kp}^{m}-D_{km...}^{ij...}L_{lp}^{m}-...
\end{equation*}%
and the vertical ($v-$) covariant derivative%
\begin{equation*}
D_{kl...}^{ij...}|_{p}=\frac{\partial D_{kl...}^{ij...}}{\partial y^{p}}%
+D_{kl...}^{mj...}C_{mp}^{i}+D_{kl...}^{im...}C_{mp}^{j}+...-D_{ml...}^{ij...}C_{kp}^{m}-D_{km...}^{ij...}C_{lp}^{m}-...,
\end{equation*}%
where%
\begin{equation*}
D=D_{kl...}^{ij...}(x,y)\frac{\delta }{\delta x^{i}}\otimes \frac{\partial }{%
\partial y^{j}}\otimes dx^{k}\otimes \delta y^{l}\otimes ...
\end{equation*}%
is an arbitrary d-tensor on $TM$.

\begin{remark}
The Cartan covariant derivatives produced by $C\Gamma (N)$ have the metrical
properties%
\begin{equation*}
g_{ij|k}=g_{\text{ \ }|k}^{ij}=0,\qquad g_{ij}|_{k}=g^{ij}|_{k}=0.
\end{equation*}
\end{remark}

For the study of relativistic magnetized non-viscous plasma in the
non-isotropic space-time $\mathcal{GL}^{n}$, one uses the following
geometrical objects \cite{Girtu-Girtu-Posto}:

\begin{enumerate}
\item the unit velocity-d-field of a test particle is given by%
\begin{equation*}
U=u^{i}(x,y)\frac{\partial }{\partial y^{i}},
\end{equation*}%
where, if we use the notation $\varepsilon ^{2}=g_{pq}y^{p}y^{q}>0$, then we
put $u^{i}=y^{i}/\varepsilon $. Obviously, we have $u_{i}u^{i}=1$, where $%
u_{i}=g_{im}u^{m}$;

\item the distinguished $2$-form of the (electric field)-(magnetic
induction) is given by%
\begin{equation*}
H=H_{ij}(x,y)dx^{i}\wedge dx^{j};
\end{equation*}

\item the distinguished $2$-form of the (electric induction)-(magnetic
field) is given by%
\begin{equation*}
G=G_{ij}(x,y)dx^{i}\wedge dx^{j};
\end{equation*}

\item the Minkowski energy d-tensor of the electromagnetic field inside the
non-isotropic plasma is given by%
\begin{equation*}
E=E_{ij}(x,y)dx^{i}\otimes dx^{j}+E_{ij}(x,y)\delta y^{i}\otimes \delta
y^{j}.
\end{equation*}%
The Minkowski energy adapted components are defined by%
\begin{equation*}
E_{ij}=\frac{1}{4}g_{ij}H_{rs}G^{rs}+g^{rs}H_{ir}G_{js},
\end{equation*}%
where $G^{rs}=g^{rp}g^{sq}G_{pq}$, and they must verify the Lorentz
conditions%
\begin{equation}
E_{i|m}^{m}u^{i}=0,\qquad E_{i}^{m}|_{m}u^{i}=0,
\label{Lorentz-conditions-Lagrange}
\end{equation}%
where $E_{i}^{m}=g^{mp}E_{pi}$. If we denote $H_{r}^{m}=g^{mp}H_{pr}$ and $%
G_{i}^{r}=g^{rs}G_{si}$, then we have%
\begin{equation*}
E_{i}^{m}=\frac{1}{4}\delta _{i}^{m}H_{rs}G^{rs}-H_{r}^{m}G_{i}^{r}.
\end{equation*}
\end{enumerate}

The energy-stress-momentum d-tensor, that characterize the relativistic
magnetized non-viscous plasma in a non-isotropic space-time, is defined by
the distinguished tensor \cite{Girtu-Girtu-Posto}%
\begin{equation*}
\mathcal{T}=\mathcal{T}_{ij}(x,y)dx^{i}\otimes dx^{j}+\mathcal{T}%
_{ij}(x,y)\delta y^{i}\otimes \delta y^{j},
\end{equation*}%
whose adapted components are%
\begin{equation}
\mathcal{T}_{ij}=\left( \mathbf{\rho }+\frac{\mathbf{p}}{c^{2}}\right)
u_{i}u_{j}+\mathbf{p}g_{ij}+E_{ij},  \label{stress-Lagrange}
\end{equation}%
where $c=$ constant, $\mathbf{p}=\mathbf{p}(x,y)$ and $\mathbf{\rho }=%
\mathbf{\rho }(x,y)$ have the similar physical meanings as in the
semi-Riemannian case.

In the Lagrangian framework of plasma, one postulates that the following 
\textit{conservation laws} for the components (\ref{stress-Lagrange}) are
true \cite{Girtu-Girtu-Posto}:%
\begin{equation}
\mathcal{T}_{i|m}^{m}=0,\qquad \mathcal{T}_{i}^{m}|_{m}=0,
\label{Cons-Lagrange}
\end{equation}%
where%
\begin{equation*}
\mathcal{T}_{i}^{m}=g^{mp}\mathcal{T}_{pi}=\left( \mathbf{\rho }+\frac{%
\mathbf{p}}{c^{2}}\right) u^{m}u_{i}+\mathbf{p}\delta _{i}^{m}+E_{i}^{m}.
\end{equation*}%
By direct computations, the conservation equations (\ref{Cons-Lagrange})
become%
\begin{equation}
\begin{array}{l}
\left[ \left( \mathbf{\rho }+\dfrac{\mathbf{p}}{c^{2}}\right) u^{m}\right]
_{|m}u_{i}+\left( \mathbf{\rho }+\dfrac{\mathbf{p}}{c^{2}}\right)
u^{m}u_{i|m}+\mathbf{p}_{,,i}-g_{ir}\overset{h}{\mathcal{F}^{r}}=0,\medskip
\\ 
\left. \left[ \left( \mathbf{\rho }+\dfrac{\mathbf{p}}{c^{2}}\right) u^{m}%
\right] \right\vert _{m}u_{i}+\left( \mathbf{\rho }+\dfrac{\mathbf{p}}{c^{2}}%
\right) u^{m}u_{i}|_{m}+\mathbf{p}_{\#i}-g_{ir}\overset{v}{\mathcal{F}^{r}}%
=0,%
\end{array}
\label{cons-local-Lagrange}
\end{equation}%
where $\mathbf{p}_{,,i}=\delta \mathbf{p}/\delta x^{i}$, $\mathbf{p}%
_{\#i}=\partial \mathbf{p}/\partial y^{i}$ and $\overset{h}{\mathcal{F}^{r}}%
=-g^{rs}E_{s|m}^{m}$ ($\overset{v}{\mathcal{F}^{r}}=-g^{rs}E_{s}^{m}|_{m}$,
respectively) is the \textit{horizontal }(\textit{vertical}, respectively) 
\textit{Lorentz force}.

Contracting the conservation equations (\ref{cons-local-Lagrange}) with $%
u^{i}$ and taking into account the Lorentz conditions (\ref%
{Lorentz-conditions-Lagrange}), we find the \textit{continuity equations} of
plasma in a non-isotropic medium:%
\begin{equation}
\begin{array}{l}
\left[ \left( \mathbf{\rho }+\dfrac{\mathbf{p}}{c^{2}}\right) u^{m}\right]
_{|m}+\mathbf{p}_{,,m}u^{m}=0,\medskip \\ 
\left. \left[ \left( \mathbf{\rho }+\dfrac{\mathbf{p}}{c^{2}}\right) u^{m}%
\right] \right\vert _{m}+\mathbf{p}_{\#m}u^{m}=0,%
\end{array}
\label{Continuity-Lagrange}
\end{equation}%
where we also used the equalities%
\begin{equation*}
\begin{array}{l}
\medskip 0=u_{i}u_{|m}^{i}=\dfrac{1}{2}\left( u_{i}u^{i}\right)
_{,,m}=-u_{i|m}u^{i}, \\ 
0=u_{i}u^{i}|_{m}=\dfrac{1}{2}\left( u_{i}u^{i}\right)
_{\#m}=-u_{i}|_{m}u^{i},%
\end{array}%
\end{equation*}%
the symbols "$_{,,m}$" and "$_{\#m}$" being the derivative operators $\delta
/\delta x^{m}$ and $\partial /\partial y^{m}$.

Replacing the continuity laws (\ref{Continuity-Lagrange}) into the
conservation equations (\ref{cons-local-Lagrange}), we find the \textit{%
relativistic Euler equations} for non-isotropic plasma, namely%
\begin{equation}
\begin{array}{l}
\medskip \left( \mathbf{\rho }+\dfrac{\mathbf{p}}{c^{2}}\right) u_{i|m}u^{m}-%
\mathbf{p}_{,,m}\left( u^{m}u_{i}-\delta _{i}^{m}\right) -g_{im}\overset{h}{%
\mathcal{F}^{m}}=0, \\ 
\left( \mathbf{\rho }+\dfrac{\mathbf{p}}{c^{2}}\right) u_{i}|_{m}u^{m}-%
\mathbf{p}_{\#m}\left( u^{m}u_{i}-\delta _{i}^{m}\right) -g_{im}\overset{v}{%
\mathcal{F}^{m}}=0.%
\end{array}
\label{Euler-Lagrange}
\end{equation}

If we take now $y^{m}=dx^{m}/dt$, then we have%
\begin{equation*}
u^{m}=\frac{1}{\varepsilon _{0}}\frac{dx^{m}}{dt}=\frac{dx^{m}}{ds},\qquad
\varepsilon _{0}^{2}=g_{ij}(x,dx/dt)\frac{dx^{i}}{dt}\frac{dx^{j}}{dt}.
\end{equation*}%
Introducing this $u^{m}$ into Euler equations (\ref{Euler-Lagrange}), we
obtain the equations of the \textit{stream lines} for non-isotropic plasma,
which are given by the following second order DE systems:

\begin{itemize}
\item \textit{horizontal} stream line DEs:%
\begin{equation*}
\begin{array}{l}
\dfrac{d^{2}x^{k}}{ds^{2}}+\left[ L_{rm}^{k}-\dfrac{c^{2}}{\mathbf{p}+%
\mathbf{\rho }c^{2}}\delta _{r}^{k}\mathbf{p}_{,,m}\right] \dfrac{dx^{r}}{ds}%
\dfrac{dx^{m}}{ds}=\dfrac{c^{2}}{\mathbf{p}+\mathbf{\rho }c^{2}}\left[ 
\overset{h}{\mathcal{F}^{k}}-g^{km}\mathbf{p}_{,,m}\right] +\medskip \\ 
\medskip +\dfrac{N_{m}^{k}}{\varepsilon _{0}}\dfrac{dx^{m}}{ds}-\dfrac{%
N_{m}^{p}g_{pr}}{\varepsilon _{0}}\dfrac{dx^{r}}{ds}\dfrac{dx^{m}}{ds}\dfrac{%
dx^{k}}{ds}- \\ 
-\dfrac{N_{m}^{r}}{2}\dfrac{\partial g_{pq}}{\partial y^{r}}\dfrac{dx^{p}}{ds%
}\dfrac{dx^{q}}{ds}\dfrac{dx^{m}}{ds}\dfrac{dx^{k}}{ds};%
\end{array}%
\end{equation*}

\item \textit{vertical} stream line DEs:%
\begin{equation*}
\begin{array}{l}
\left[ C_{rm}^{k}-\dfrac{c^{2}}{\mathbf{p}+\mathbf{\rho }c^{2}}\delta
_{r}^{k}\mathbf{p}_{\#m}\right] \dfrac{dx^{r}}{ds}\dfrac{dx^{m}}{ds}=\dfrac{%
c^{2}}{\mathbf{p}+\mathbf{\rho }c^{2}}\left[ \overset{v}{\mathcal{F}^{k}}%
-g^{km}\mathbf{p}_{\#m}\right] +\medskip \\ 
+\dfrac{1}{2}\dfrac{\partial g_{pq}}{\partial y^{r}}\dfrac{dx^{p}}{ds}\dfrac{%
dx^{q}}{ds}\dfrac{dx^{r}}{ds}\dfrac{dx^{k}}{ds}.%
\end{array}%
\end{equation*}
\end{itemize}

\begin{remark}
If the the metrical d-tensor $g_{ij}(x,y)$ is Finslerian one, that is we have%
\begin{equation*}
g_{ij}(x,y)=\frac{1}{2}\frac{\partial ^{2}F^{2}}{\partial y^{i}\partial y^{j}%
},
\end{equation*}%
where $F:TM\rightarrow \mathbb{R}_{+}$ is a Finslerian metric, then the DEs
of \textbf{stream lines} of plasma in non-isotropic spaces reduce to

\begin{itemize}
\item \textbf{horizontal} stream line DEs:%
\begin{equation*}
\begin{array}{l}
\dfrac{d^{2}x^{k}}{ds^{2}}+\left[ L_{rm}^{k}-\dfrac{c^{2}}{\mathbf{p}+%
\mathbf{\rho }c^{2}}\delta _{r}^{k}\mathbf{p}_{,,m}\right] \dfrac{dx^{r}}{ds}%
\dfrac{dx^{m}}{ds}=\dfrac{c^{2}}{\mathbf{p}+\mathbf{\rho }c^{2}}\left[ 
\overset{h}{\mathcal{F}^{k}}-g^{km}\mathbf{p}_{,,m}\right] +\medskip \\ 
+\dfrac{2}{F^{2}}\left[ G^{k}-g_{pr}G^{p}\dfrac{dx^{r}}{ds}\dfrac{dx^{k}}{ds}%
\right] ;%
\end{array}%
\end{equation*}

\item \textbf{vertical} stream line DEs:%
\begin{equation*}
\mathbf{p}_{\#m}\left[ g^{mk}-\dfrac{dx^{m}}{ds}\dfrac{dx^{k}}{ds}\right] =%
\overset{v}{\mathcal{F}^{k}},
\end{equation*}%
where, if the generalized Christoffel symbols of $g_{ij}(x,y)$ are%
\begin{equation*}
\Gamma _{jk}^{i}(x,y)=\frac{g^{im}}{2}\left( \frac{\partial g_{jm}}{\partial
x^{k}}+\frac{\partial g_{km}}{\partial x^{j}}-\frac{\partial g_{jk}}{%
\partial x^{m}}\right) ,
\end{equation*}%
then we have%
\begin{equation*}
G^{k}=\frac{1}{2}\Gamma _{pq}^{i}(x,y)y^{p}y^{q}.
\end{equation*}
\end{itemize}
\end{remark}

\section{The Riemann-Lagrange geometrical approach. Multi-time plasma}

\hspace{5mm}Let us consider that $(T^{p},h_{\alpha \beta }(t))$ is a
Riemannian manifold of dimension $p$, whose local coordinates are $%
(t^{\alpha })_{\alpha =\overline{1,p}}$. Suppose that the Christoffel
symbols of the Riemannian metric $h_{\alpha \beta }(t)$ are $\varkappa
_{\alpha \beta }^{\gamma }(t)$. Let $J^{1}(T,M)$ be the $1$-jet space (it
has the dimension $p+n+pn$) whose local coordinates are $(t^{\alpha
},x^{i},x_{\alpha }^{i})$. These transform by the rules 
\begin{equation*}
\left\{ 
\begin{array}{l}
\widetilde{t}^{\alpha }=\widetilde{t}^{\alpha }(t^{\beta })\medskip \\ 
\widetilde{x}^{i}=\widetilde{x}^{i}(x^{j})\medskip \\ 
\widetilde{x}_{\alpha }^{i}={{\dfrac{\partial \widetilde{x}^{i}}{\partial
x^{j}}}{\dfrac{\partial t^{\beta }}{\partial \widetilde{t}^{\alpha }}}%
x_{\beta }^{j},}%
\end{array}%
\right.
\end{equation*}%
where $\det (\partial \widetilde{t}^{\alpha }/\partial t^{\beta })\neq 0$
and $\det (\partial \widetilde{x}^{i}/\partial x^{j})\neq 0$. Note that,
throughout this work, the greek letters run from $1$ to $p$ and the latin
letters run from $1$ to $n$.

Let $\mathcal{GML}_{p}^{n}=\left( J^{1}(T,M),\text{ }G_{(i)(j)}^{(\alpha
)(\beta )}=h^{\alpha \beta }g_{ij}\right) $ be a multi-time generalized
Lagrange space (for more details, please see Neagu \cite{Neagu-Forum}),
where $g_{ij}(t^{\gamma },x^{k},x_{\gamma }^{k})$ is a metrical d-tensor on $%
J^{1}(T,M)$, which is symmetrical, non-degenerate and has a constant
signature.

Let us consider that $\mathcal{GML}_{p}^{n}$ is endowed with a nonlinear
connection having the form \cite{Neagu-Forum}%
\begin{equation*}
\Gamma =\left( M_{(\alpha )\beta }^{(i)}=-\varkappa _{\alpha \beta }^{\gamma
}x_{\gamma }^{i},\text{ }N_{(\alpha )j}^{(i)}\right) .
\end{equation*}%
The nonlinear connection $\Gamma $ produces on $J^{1}(T,M)$ the following
dual adapted bases of d-vectors and d-covectors:%
\begin{equation*}
\left\{ \frac{\delta }{\delta t^{\alpha }},\frac{\delta }{\delta x^{i}},%
\frac{\partial }{\partial x_{\alpha }^{i}}\right\} \subset \mathcal{X}%
(J^{1}(T,M)),\qquad \left\{ dt^{\alpha },dx^{i},\delta x_{\alpha
}^{i}\right\} \subset \mathcal{X}^{\ast }(J^{1}(T,M)),
\end{equation*}%
where%
\begin{equation*}
\frac{\delta }{\delta t^{\alpha }}=\frac{\partial }{\partial t^{\alpha }}%
+\varkappa _{\alpha \mu }^{\gamma }x_{\gamma }^{m}\frac{\partial }{\partial
x_{\mu }^{m}},\qquad \frac{\delta }{\delta x^{i}}=\frac{\partial }{\partial
x^{i}}-N_{(\mu )i}^{(m)}\frac{\partial }{\partial x_{\mu }^{m}},
\end{equation*}%
\begin{equation*}
\delta x_{\alpha }^{i}=dx_{\alpha }^{i}-\varkappa _{\alpha \mu }^{\gamma
}x_{\gamma }^{i}dt^{\mu }+N_{(\alpha )m}^{(i)}dx^{m}.
\end{equation*}%
Note that the d-tensors on the $1$-jet space $J^{1}(T,M)$ also behave like
classical tensors. For example, on the $1$-jet space $J^{1}(T,M)$ we have
the global metrical d-tensor%
\begin{equation*}
\mathbb{G}=h_{\alpha \beta }dt^{\alpha }\otimes dt^{\beta
}+g_{ij}dx^{i}\otimes dx^{j}+h^{\alpha \beta }g_{ij}\delta x_{\alpha
}^{i}\otimes \delta x_{\beta }^{j},
\end{equation*}%
which may be endowed with the physical meaning of non-isotropic multi-time
gravitational potential. It follows that $\mathbb{G}$ has the adapted
components%
\begin{equation*}
\mathbb{G}_{AB}=\left\{ 
\begin{array}{llll}
h_{\alpha \beta }, & \text{for} & A=\alpha , & B=\beta \medskip \\ 
g_{ij}, & \text{for} & A=i, & B=j\medskip \\ 
h^{\alpha \beta }g_{ij}, & \text{for} & A=_{(i)}^{(\alpha )}, & 
B=_{(j)}^{(\beta )}\medskip \\ 
0, & \text{otherwise.} &  & 
\end{array}%
\right.
\end{equation*}

Following the geometrical ideas of Asanov \cite{Asanov} and Neagu \cite%
{Neagu-Forum}, the preceding geometrical ingredients lead us to the the
Cartan canonical $\Gamma $-linear connection%
\begin{equation*}
C\Gamma =\left( \varkappa _{\alpha \beta }^{\gamma },\text{ }G_{j\gamma
}^{k},\text{ }L_{jk}^{i},\text{ }C_{j(k)}^{i(\gamma )}\right) ,
\end{equation*}%
where%
\begin{equation}
\begin{array}{c}
\medskip G_{j\gamma }^{k}=\dfrac{g^{km}}{2}\dfrac{\delta g_{mj}}{\delta
t^{\gamma }},\qquad L_{jk}^{i}=\dfrac{g^{im}}{2}\left( \dfrac{\delta g_{jm}}{%
\delta x^{k}}+\dfrac{\delta g_{km}}{\delta x^{j}}-\dfrac{\delta g_{jk}}{%
\delta x^{m}}\right) , \\ 
C_{j(k)}^{i(\gamma )}=\dfrac{g^{im}}{2}\left( \dfrac{\partial g_{jm}}{%
\partial x_{\gamma }^{k}}+\dfrac{\partial g_{km}}{\partial x_{\gamma }^{j}}-%
\dfrac{\partial g_{jk}}{\partial x_{\gamma }^{m}}\right) .%
\end{array}
\label{Cartan-Multi-Time}
\end{equation}%
In the sequel, the Cartan linear connection $C\Gamma $, given by (\ref%
{Cartan-Multi-Time}), induces the $T$-horizontal ($h_{T}-$) covariant
derivative%
\begin{eqnarray*}
D_{\gamma k(\beta )(l).../\varepsilon }^{\alpha i(j)(\nu )...} &=&\frac{%
\delta D_{\gamma k(\beta )(l)...}^{\alpha i(j)(\nu )...}}{\delta
t^{\varepsilon }}+D_{\gamma k(\beta )(l)...}^{\mu i(j)(\nu )...}\varkappa
_{\mu \varepsilon }^{\alpha }+D_{\gamma k(\beta )(l)...}^{\alpha m(j)(\nu
)...}G_{m\varepsilon }^{i}+ \\
&&+D_{\gamma k(\beta )(l)...}^{\alpha i(m)(\nu )...}G_{m\varepsilon
}^{j}+D_{\gamma k(\beta )(l)...}^{\alpha i(j)(\mu )...}\varkappa _{\mu
\varepsilon }^{\nu }+...- \\
&&-D_{\mu k(\beta )(l)...}^{\alpha i(j)(\nu )...}\varkappa _{\gamma
\varepsilon }^{\mu }-D_{\gamma m(\beta )(l)...}^{\alpha i(j)(\nu
)...}G_{k\varepsilon }^{m}- \\
&&-D_{\gamma k(\mu )(l)...}^{\alpha i(j)(\nu )...}\varkappa _{\beta
\varepsilon }^{\mu }-D_{\gamma k(\beta )(m)...}^{\alpha i(j)(\nu
)...}G_{l\varepsilon }^{m}...,
\end{eqnarray*}%
the $M$-horizontal ($h_{M}-$) covariant derivative%
\begin{eqnarray*}
D_{\gamma k(\beta )(l)...|p}^{\alpha i(j)(\nu )...} &=&\frac{\delta
D_{\gamma k(\beta )(l)...}^{\alpha i(j)(\nu )...}}{\delta x^{p}}+D_{\gamma
k(\beta )(l)...}^{\alpha m(j)(\nu )...}L_{mp}^{i}+D_{\gamma k(\beta
)(l)...}^{\alpha i(m)(\nu )...}L_{mp}^{j}+...- \\
&&-D_{\gamma m(\beta )(l)...}^{\alpha i(j)(\nu )...}L_{kp}^{m}-D_{\gamma
k(\beta )(m)...}^{\alpha i(j)(\nu )...}L_{lp}^{m}-...
\end{eqnarray*}%
and the vertical ($v-$) covariant derivative%
\begin{eqnarray*}
D_{\gamma k(\beta )(l)...}^{\alpha i(j)(\nu )...}|_{(p)}^{(\varepsilon )} &=&%
\frac{\partial D_{\gamma k(\beta )(l)...}^{\alpha i(j)(\nu )...}}{\partial
x_{\varepsilon }^{p}}+D_{\gamma k(\beta )(l)...}^{\alpha m(j)(\nu
)...}C_{m(p)}^{i(\varepsilon )}+D_{\gamma k(\beta )(l)...}^{\alpha i(m)(\nu
)...}C_{m(p)}^{j(\varepsilon )}+...- \\
&&-D_{\gamma m(\beta )(l)...}^{\alpha i(j)(\nu )...}C_{k(p)}^{m(\varepsilon
)}-D_{\gamma k(\beta )(m)...}^{\alpha i(j)(\nu )...}C_{l(p)}^{m(\varepsilon
)}-...,
\end{eqnarray*}%
where%
\begin{equation*}
D=D_{\gamma k(\beta )(l)...}^{\alpha i(j)(\nu )...}(t^{\lambda
},x^{r},x_{\lambda }^{r})\frac{\delta }{\delta t^{\alpha }}\otimes \frac{%
\delta }{\delta x^{i}}\otimes \frac{\partial }{\partial x_{\beta }^{j}}%
\otimes dt^{\gamma }\otimes dx^{k}\otimes \delta x_{\nu }^{l}\otimes ...
\end{equation*}%
is an arbitrary d-tensor on $J^{1}(T,M)$.

\begin{remark}
The Cartan covariant derivatives produced by $C\Gamma $ have the
me\-tri\-cal properties%
\begin{equation*}
\begin{array}{lll}
h_{\alpha \beta /\gamma }=h_{\text{ \ }/\gamma }^{\alpha \beta }=0, & 
h_{\alpha \beta |k}=h_{\text{ \ }|k}^{\alpha \beta }=0, & h_{\alpha \beta
}|_{(k)}^{(\gamma )}=h^{\alpha \beta }|_{(k)}^{(\gamma )}=0,\medskip \\ 
g_{ij/\gamma }=g_{\text{ \ }/\gamma }^{ij}=0, & g_{ij|k}=g_{\text{ \ }%
|k}^{ij}=0, & g_{ij}|_{(k)}^{(\gamma )}=g^{ij}|_{(k)}^{(\gamma )}=0.%
\end{array}%
\end{equation*}
\end{remark}

For the study of the relativistic magnetized non-viscous plasma dynamics, in
a Riemann-Lagrange geometrical multi-time approach, we use the following
geometrical objects:

\begin{enumerate}
\item the unit multi-time velocity-d-field of a test particle is given by%
\begin{equation*}
U=u_{\alpha }^{i}(t^{\gamma },x^{k},x_{\gamma }^{k})\frac{\partial }{%
\partial x_{\alpha }^{i}},
\end{equation*}%
where, if we take $\varepsilon ^{2}=h^{\mu \nu }g_{pq}x_{\mu }^{p}x_{\nu
}^{q}>0$, then we put $u_{\alpha }^{i}=x_{\alpha }^{i}/\varepsilon $.
Obviously, we have $h^{\alpha \beta }u_{i\alpha }u_{\beta }^{i}=1$, where $%
u_{i\alpha }=g_{im}u_{\alpha }^{m}$;

\item the distinguished multi-time $2$-form of the (electric
field)-(magnetic induction) is given by%
\begin{equation*}
H=H_{ij}(t^{\gamma },x^{k},x_{\gamma }^{k})dx^{i}\wedge dx^{j};
\end{equation*}

\item the distinguished multi-time $2$-form of the (electric
induction)-(magnetic field) is given by%
\begin{equation*}
G=G_{ij}(t^{\gamma },x^{k},x_{\gamma }^{k})dx^{i}\wedge dx^{j};
\end{equation*}

\item the multi-time Minkowski energy d-tensor of the electromagnetic field
inside the multi-time plasma is given by%
\begin{equation*}
E=E_{ij}(t^{\gamma },x^{k},x_{\gamma }^{k})dx^{i}\otimes dx^{j}+h^{\eta \nu
}E_{ij}(t^{\gamma },x^{k},x_{\gamma }^{k})\delta x_{\eta }^{i}\otimes \delta
x_{\nu }^{j}.
\end{equation*}%
The multi-time Minkowski energy adapted components are defined by the
similar formulas%
\begin{equation*}
E_{ij}=\dfrac{1}{4}g_{ij}H_{rs}G^{rs}+g^{rs}H_{ir}G_{js},
\end{equation*}%
where $G^{rs}=g^{rp}g^{sq}G_{pq}$. Furthermore, we suppose that the
multi-time Min\-kow\-ski energy adapted components verify the multi-time
Lorentz con\-di\-tions%
\begin{equation}
E_{i|m}^{m}u_{\alpha }^{i}=0,\qquad E_{i}^{m}|_{(m)}^{(\mu )}u_{\mu }^{i}=0,
\label{Lorentz-conditions-Multi-Time}
\end{equation}%
where $E_{i}^{m}=g^{mp}E_{pi}$. Obviously, if we use the notations $%
H_{r}^{m}=g^{mp}H_{pr}$ and $G_{i}^{r}=g^{rs}G_{si}$, we obtain%
\begin{equation*}
E_{i}^{m}=\frac{1}{4}\delta _{i}^{m}H_{rs}G^{rs}-H_{r}^{m}G_{i}^{r}.
\end{equation*}
\end{enumerate}

In our Riemann-Lagrange geometrical approach, the multi-time plasma is
characterized by the energy-stress-momentum d-tensor defined by%
\begin{equation*}
\mathcal{T}=\mathcal{T}_{ij}(t^{\gamma },x^{k},x_{\gamma }^{k})dx^{i}\otimes
dx^{j}+h^{\eta \nu }\mathcal{T}_{ij}(t^{\gamma },x^{k},x_{\gamma
}^{k})\delta x_{\eta }^{i}\otimes \delta x_{\nu }^{j},
\end{equation*}%
where%
\begin{equation}
\mathcal{T}_{ij}=\left( \mathbf{\rho }+\dfrac{\mathbf{p}}{c^{2}}\right)
h^{\alpha \beta }u_{i\alpha }u_{j\beta }+\mathbf{p}g_{ij}+E_{ij}.
\label{stress-Multi-Time-1}
\end{equation}%
The entities $c=$ constant, $\mathbf{p}=\mathbf{p}(t^{\gamma
},x^{k},x_{\gamma }^{k})$ and $\mathbf{\rho }=\mathbf{\rho }(t^{\gamma
},x^{k},x_{\gamma }^{k})$ have the multi-time extended physical meanings of
their analogous entities from the semi-Riemannian framework. Note that the
adapted components of the energy-stress-momentum d-tensor $\mathcal{T}$ of
multi-time plasma are given by%
\begin{equation}
\mathcal{T}_{CF}=\left\{ 
\begin{array}{llll}
\mathcal{T}_{ij}, & \text{for} & C=i, & F=j\medskip \\ 
h^{\eta \nu }\mathcal{T}_{ij}, & \text{for} & C=_{(i)}^{(\eta )}, & 
F=_{(j)}^{(\nu )}\medskip \\ 
0, & \text{otherwise.} &  & 
\end{array}%
\right.  \label{stress-Multi-Time-2}
\end{equation}

In the multi-time Riemann-Lagrange framework of plasma, we postulate that
the following \textit{multi-time conservation laws} for the components (\ref%
{stress-Multi-Time-1}) and (\ref{stress-Multi-Time-2}) are true:%
\begin{equation}
\mathcal{T}_{A:M}^{M}=0,\qquad \forall \text{ }A\in \left\{ \alpha ,\text{ }%
i,\text{ }_{(i)}^{(\alpha )}\right\} ,  \label{cons-general-M-T}
\end{equation}%
where the capital latin letters $A,M,...$ are indices of kind $\alpha ,$ $i$
or $_{(i)}^{(\alpha )}$, "$_{:M}$" represents one of the local covariant
derivatives $h_{T}-$, $h_{M}-$ or $v-$ and%
\begin{equation*}
\mathcal{T}_{A}^{M}=\mathbb{G}^{MD}\mathcal{T}_{DA}=\left\{ 
\begin{array}{llll}
\mathcal{T}_{i}^{m}, & \text{for} & A=i, & M=m\medskip \\ 
\delta _{\mu }^{\alpha }\mathcal{T}_{i}^{m}, & \text{for} & 
A=_{(i)}^{(\alpha )}, & M=_{(m)}^{(\mu )}\medskip \\ 
0, & \text{otherwise.} &  & 
\end{array}%
\right.
\end{equation*}%
Obviously, the d-tensor $\mathcal{T}_{i}^{m}$ is given by the formula%
\begin{equation*}
\mathcal{T}_{i}^{m}=g^{mp}\mathcal{T}_{pi}=\left( \mathbf{\rho }+\dfrac{%
\mathbf{p}}{c^{2}}\right) h^{\alpha \beta }u_{\alpha }^{m}u_{i\beta }+%
\mathbf{p}\delta _{i}^{m}+E_{i}^{m}.
\end{equation*}%
The multi-time conservation laws (\ref{cons-general-M-T}) reduce to the
multi-time conservation equations%
\begin{equation}
\mathcal{T}_{i|m}^{m}=0,\qquad \mathcal{T}_{i}^{m}|_{(m)}^{(\mu )}=0.
\label{cons-eq-M-T}
\end{equation}%
By direct computations, the multi-time conservation equations (\ref%
{cons-eq-M-T}) become%
\begin{equation}
\begin{array}{l}
h^{\alpha \beta }\left[ \left( \mathbf{\rho }+\dfrac{\mathbf{p}}{c^{2}}%
\right) u_{\alpha }^{m}\right] _{|m}u_{i\beta }+\left( \mathbf{\rho }+\dfrac{%
\mathbf{p}}{c^{2}}\right) h^{\alpha \beta }u_{\alpha }^{m}u_{i\beta |m}+%
\mathbf{p}_{,,i}-g_{ir}\overset{h}{\mathcal{F}^{r}}=0,\medskip \\ 
h^{\alpha \beta }\left. \left[ \left( \mathbf{\rho }+\dfrac{\mathbf{p}}{c^{2}%
}\right) u_{\alpha }^{m}\right] \right\vert _{(m)}^{(\mu )}u_{i\beta
}+\left( \mathbf{\rho }+\dfrac{\mathbf{p}}{c^{2}}\right) h^{\alpha \beta
}u_{\alpha }^{m}u_{i\beta }|_{(m)}^{(\mu )}+\medskip \\ 
+\mathbf{p}_{\#(i)}^{\text{ \ }(\mu )}-g_{ir}\overset{v}{\mathcal{F}}\text{ }%
\!\!^{r\mu }=0,%
\end{array}
\label{cons-local-M-T}
\end{equation}%
where $\mathbf{p}_{,,i}=\delta \mathbf{p}/\delta x^{i}$, $\mathbf{p}%
_{\#(i)}^{\text{ \ }(\mu )}=\partial \mathbf{p}/\partial x_{\mu }^{i}$ and

\begin{itemize}
\item $\overset{h}{\mathcal{F}^{r}}=-g^{rs}E_{s|m}^{m}$ is the \textit{%
multi-time horizontal Lorentz force};

\item $\overset{v}{\mathcal{F}}$ $\!\!^{r\mu }=-g^{rs}E_{s}^{m}|_{(m)}^{(\mu
)}$ is the \textit{multi-time vertical} \textit{Lorentz d-tensor force.}
\end{itemize}

Contracting the multi-time conservation equations (\ref{cons-local-M-T})
with $u_{\mu }^{i}$ and taking into account the Lorentz conditions (\ref%
{Lorentz-conditions-Multi-Time}), we find the \textit{continuity equations}
of multi-time plasma, namely%
\begin{equation*}
\begin{array}{l}
h^{\alpha \beta }\left[ \left( \mathbf{\rho }+\dfrac{\mathbf{p}}{c^{2}}%
\right) u_{\alpha }^{m}\right] _{|m}u_{i\beta }u_{\mu }^{i}+\left( \mathbf{%
\rho }+\dfrac{\mathbf{p}}{c^{2}}\right) h^{\alpha \beta }u_{\alpha
}^{m}u_{i\beta |m}u_{\mu }^{i}+\mathbf{p}_{,,m}u_{\mu }^{m}=0,\medskip \\ 
h^{\alpha \beta }\left. \left[ \left( \mathbf{\rho }+\dfrac{\mathbf{p}}{c^{2}%
}\right) u_{\alpha }^{m}\right] \right\vert _{(m)}^{(\mu )}u_{i\beta }u_{\mu
}^{i}+\left( \mathbf{\rho }+\dfrac{\mathbf{p}}{c^{2}}\right) h^{\alpha \beta
}u_{\alpha }^{m}u_{i\beta }|_{(m)}^{(\mu )}u_{\mu }^{i}+\mathbf{p}_{\#(m)}^{%
\text{ \ }(\mu )}u_{\mu }^{m}=0.%
\end{array}%
\end{equation*}

If we take now $x_{\eta }^{l}=\partial x^{l}/dt^{\eta }$, then we have%
\begin{equation*}
u_{\eta }^{l}=\frac{x_{\eta }^{l}}{\varepsilon _{0}},\qquad \varepsilon
_{0}^{2}=h^{\alpha \beta }(t)g_{ij}(t^{\gamma },x^{k},x_{\gamma
}^{k})x_{\alpha }^{i}x_{\beta }^{j}.
\end{equation*}%
Introducing this $u_{\eta }^{l}$ into multi-time conservation equations (\ref%
{cons-local-M-T}), we obtain the second order PDEs of the \textit{stream
sheets} that characterize the multi-time plasma:

\begin{itemize}
\item \textit{horizontal} stream sheet PDEs:%
\begin{equation*}
\begin{array}{l}
h^{\alpha \beta }\left\{ \left[ \left( \mathbf{\rho }+\dfrac{\mathbf{p}}{%
c^{2}}\right) \dfrac{x_{\alpha }^{m}}{\varepsilon _{0}}\right] _{|m}x_{\beta
}^{k}+\left( \mathbf{\rho }+\dfrac{\mathbf{p}}{c^{2}}\right) x_{\alpha }^{m}%
\left[ \dfrac{x_{\beta }^{k}}{\varepsilon _{0}}\right] _{|m}\right\}
=\medskip \\ 
=\varepsilon _{0}\left[ \overset{h}{\mathcal{F}^{k}}-g^{km}\mathbf{p}_{,,m}%
\right] ;%
\end{array}%
\end{equation*}

\item \textit{vertical} stream sheet PDEs:%
\begin{equation*}
\begin{array}{l}
h^{\alpha \beta }\left\{ \left. \left[ \left( \mathbf{\rho }+\dfrac{\mathbf{p%
}}{c^{2}}\right) \dfrac{x_{\alpha }^{m}}{\varepsilon _{0}}\right]
\right\vert _{(m)}^{(\mu )}x_{\beta }^{k}+\left( \mathbf{\rho }+\dfrac{%
\mathbf{p}}{c^{2}}\right) x_{\alpha }^{m}\left. \left[ \dfrac{x_{\beta }^{k}%
}{\varepsilon _{0}}\right] \right\vert _{(m)}^{(\mu )}\right\} =\medskip \\ 
=\varepsilon _{0}\left[ \overset{v}{\mathcal{F}}\text{ }\!\!^{k\mu }-g^{km}%
\mathbf{p}_{\#(m)}^{\text{ \ }(\mu )}\right] .%
\end{array}%
\end{equation*}
\end{itemize}

Taking into account the local form of the $h_{M}-$ and $v-$ covariant
derivatives produced by the Cartan connection $C\Gamma $, the expressions of
the PDEs of the stream sheets of multi-time plasma reduce to:

\begin{itemize}
\item \textit{horizontal} stream sheet PDEs:%
\begin{equation*}
\begin{array}{l}
h^{\alpha \beta }\left\{ \mathcal{H}_{m}x_{\alpha }^{m}x_{\beta }^{k}+\dfrac{%
1}{\varepsilon _{0}}\left( \mathbf{\rho }+\dfrac{\mathbf{p}}{c^{2}}\right) %
\left[ L_{rm}^{k}x_{\beta }^{r}-N_{(\beta )m}^{(k)}\right] x_{\alpha
}^{m}+\right. \medskip \\ 
\left. +\dfrac{1}{\varepsilon _{0}}\left( \mathbf{\rho }+\dfrac{\mathbf{p}}{%
c^{2}}\right) \left[ L_{rm}^{m}x_{\alpha }^{r}-N_{(\alpha )m}^{(m)}\right]
x_{\beta }^{k}\right\} =\varepsilon _{0}\left[ \overset{h}{\mathcal{F}^{k}}%
-g^{km}\mathbf{p}_{,,m}\right] ,%
\end{array}%
\end{equation*}%
where%
\begin{equation*}
\mathcal{H}_{m}=\left[ \dfrac{1}{\varepsilon _{0}}\left( \mathbf{\rho }+%
\dfrac{\mathbf{p}}{c^{2}}\right) \right] _{,,m}+\left( \mathbf{\rho }+\dfrac{%
\mathbf{p}}{c^{2}}\right) \left[ \dfrac{1}{\varepsilon _{0}}\right] _{,,m};
\end{equation*}

\item \textit{vertical} stream sheet PDEs:%
\begin{equation*}
\begin{array}{l}
h^{\alpha \beta }\left\{ \mathcal{V}_{(m)}^{(\mu )}x_{\alpha }^{m}x_{\beta
}^{k}+\dfrac{1}{\varepsilon _{0}}\left( \mathbf{\rho }+\dfrac{\mathbf{p}}{%
c^{2}}\right) \left[ n\cdot \delta _{\alpha }^{\mu }x_{\beta }^{k}+\delta
_{\beta }^{\mu }x_{\alpha }^{k}\right] +\right. \medskip \\ 
\left. +\dfrac{1}{\varepsilon _{0}}\left( \mathbf{\rho }+\dfrac{\mathbf{p}}{%
c^{2}}\right) \left[ C_{m(r)}^{k(\mu )}x_{\beta }^{m}+C_{r(m)}^{m(\mu
)}x_{\beta }^{k}\right] x_{\alpha }^{r}\right\} =\varepsilon _{0}\left[ 
\overset{v}{\mathcal{F}}\text{ }\!\!^{k\mu }-g^{km}\mathbf{p}_{\#(m)}^{\text{
\ }(\mu )}\right] ,%
\end{array}%
\end{equation*}%
where $n=\dim M$ and%
\begin{equation*}
\mathcal{V}_{(m)}^{(\mu )}=\left[ \dfrac{1}{\varepsilon _{0}}\left( \mathbf{%
\rho }+\dfrac{\mathbf{p}}{c^{2}}\right) \right] _{\#(m)}^{\text{ \ \ }(\mu
)}+\left( \mathbf{\rho }+\dfrac{\mathbf{p}}{c^{2}}\right) \left[ \dfrac{1}{%
\varepsilon _{0}}\right] _{\#(m)}^{\text{ \ \ }(\mu )}.
\end{equation*}
\end{itemize}

\section{Conclusion}

\hspace{5mm}The Riemann-Lagrange geometrical theory upon the Multi-Time
Plasma Physics may be applied for a lot of interesting multi-time
generalized Lagrange spaces with physical connotations \cite{Neagu-Forum}:

\subsection{The geometrical model $\mathcal{GRGML}_{p}^{n}$ for multi-time
General Re\-lativity and Electromagnetism}

\hspace{5mm}This generalized multi-time Lagrange space is characterized by
the fundamental metrical d-tensor%
\begin{equation*}
G_{(i)(j)}^{(\alpha )(\beta )}=h^{\alpha \beta }(t^{\gamma })e^{2\sigma
\left( t^{\gamma },x^{k},x_{\gamma }^{k}\right) }\varphi _{ij}\left(
x^{k}\right)
\end{equation*}%
and the nonlinear connection%
\begin{equation*}
\mathring{\Gamma}=\left( M_{(\alpha )\beta }^{(i)}=-\varkappa _{\alpha \beta
}^{\mu }x_{\mu }^{i},\text{ }N_{(\alpha )j}^{(i)}=\gamma _{jm}^{i}x_{\alpha
}^{m}\right) .
\end{equation*}

\subsection{The geometrical model $\mathcal{RGOGML}_{p}^{n}$ for multi-time
Relativistic Optics}

\hspace{5mm}This generalized multi-time Lagrange space is characterized by
the fundamental metrical d-tensor%
\begin{equation*}
G_{(i)(j)}^{(\alpha )(\beta )}=h^{\alpha \beta }(t^{\gamma })\left\{ \varphi
_{ij}\left( x^{k}\right) +\left[ 1-\dfrac{1}{n(t^{\gamma },x^{k},x_{\gamma
}^{k})}\right] Y_{i}Y_{j}\right\}
\end{equation*}%
and, again, by the nonlinear connection $\mathring{\Gamma}$, where $%
Y_{i}=\varphi _{im}\left( x^{k}\right) x_{\mu }^{m}X^{\mu }(t^{\gamma }).$

\subsection{The geometrical model $\mathcal{EDML}_{p}^{n}$\ for multi-time
E\-lec\-tro\-dy\-na\-mics}

\hspace{5mm}This multi-time Lagrange space is characterized by the
Lagrangian function (for more details, please see \cite{Neagu-EDML})%
\begin{equation*}
L_{ED}=h^{\alpha \beta }(t^{\gamma })\varphi _{ij}\left( x^{k}\right)
x_{\alpha }^{i}x_{\beta }^{j}+U_{(i)}^{(\alpha )}\left( t^{\gamma
},x^{k}\right) x_{\alpha }^{i}+\Phi \left( t^{\gamma },x^{k}\right) ,
\end{equation*}%
which produces the fundamental metrical d-tensor 
\begin{equation*}
G_{(i)(j)}^{(\alpha )(\beta )}=\frac{1}{2}\frac{\partial ^{2}L_{ED}}{%
\partial x_{\alpha }^{i}\partial x_{\beta }^{j}}=h^{\alpha \beta }\varphi
_{ij}
\end{equation*}%
and the nonlinear connection whose components are $M_{(\alpha )\beta
}^{(i)}=-\varkappa _{\alpha \beta }^{\mu }x_{\mu }^{i}$ and 
\begin{equation*}
N_{(\alpha )j}^{(i)}=\gamma _{jm}^{i}x_{\alpha }^{m}+\frac{h_{\alpha \mu
}\varphi ^{im}}{4}\left[ \frac{\partial U_{(m)}^{(\mu )}}{\partial x^{j}}-%
\frac{\partial U_{(j)}^{(\mu )}}{\partial x^{m}}\right] .
\end{equation*}

\subsection{The geometrical model $\mathcal{BSML}_{p}^{n}$\ for Bosonic
Strings}

\hspace{5mm}This is the multi-time Lagrange space corresponding to the
multi-time Lagrangian function%
\begin{equation*}
L_{BS}=h^{\alpha \beta }(t^{\gamma })\varphi _{ij}\left( x^{k}\right)
x_{\alpha }^{i}x_{\beta }^{j}.
\end{equation*}%
In this particular case, we have the the fundamental metrical d-tensor%
\begin{equation*}
G_{(i)(j)}^{(\alpha )(\beta )}=\frac{1}{2}\frac{\partial ^{2}L_{BS}}{%
\partial x_{\alpha }^{i}\partial x_{\beta }^{j}}=h^{\alpha \beta }\varphi
_{ij}
\end{equation*}%
and the canonical nonlinear connection $\mathring{\Gamma}$. Moreover, the
Cartan canonical connection has the following simple form:%
\begin{equation*}
C\mathring{\Gamma}=\left( \varkappa _{\alpha \beta }^{\gamma },\text{ }0,%
\text{ }\gamma _{jk}^{i},\text{ }0\right) .
\end{equation*}%
It follows that, for the multi-time Lagrange space $\mathcal{BSML}_{p}^{n}$,
the PDEs of the \textit{stream sheets} for multi-time plasma simplify as
follows:

\begin{itemize}
\item \textit{horizontal} stream sheet PDEs:%
\begin{equation*}
h^{\alpha \beta }\mathcal{H}_{m}x_{\alpha }^{m}x_{\beta }^{k}=\varepsilon
_{0}\left[ \overset{h}{\mathcal{F}^{k}}-g^{km}\mathbf{p}_{,,m}\right] ;
\end{equation*}

\item \textit{vertical} stream sheet PDEs:%
\begin{equation*}
\begin{array}{l}
h^{\alpha \beta }\left\{ \mathcal{V}_{(m)}^{(\mu )}x_{\alpha }^{m}x_{\beta
}^{k}+\dfrac{1}{\varepsilon _{0}}\left( \mathbf{\rho }+\dfrac{\mathbf{p}}{%
c^{2}}\right) \left[ n\cdot \delta _{\alpha }^{\mu }x_{\beta }^{k}+\delta
_{\beta }^{\mu }x_{\alpha }^{k}\right] \right\} =\medskip \\ 
=\varepsilon _{0}\left[ \overset{v}{\mathcal{F}}\text{ }\!\!^{k\mu }-g^{km}%
\mathbf{p}_{\#(m)}^{\text{ \ }(\mu )}\right] .%
\end{array}%
\end{equation*}
\end{itemize}

\textbf{Open Problem. }There exist real physical interpretations for our
multi-time Riemann-Lagrange geometric dynamics of plasma ?

\textbf{Author's address:} Mircea N{\scriptsize EAGU}

University Transilvania of Bra\c{s}ov, Faculty of Mathematics and Informatics

Department of Algebra, Geometry and Differential Equations

B-dul Iuliu Maniu, Nr. 50, BV 500091, Bra\c{s}ov, Romania.

\textbf{E-mails:} mircea.neagu@unitbv.ro, mirceaneagu73@yahoo.com

\textbf{Website:} http://www.2collab.com/user:mirceaneagu

\end{document}